\documentclass[11pt]{article}
\renewcommand{\baselinestretch}{1.3}

\makeatletter
\usepackage{psfig,graphicx}
\makeatletter

\def\singlespace{\def\baselinestretch{1}\@normalsize}

\newcommand{\fig}[3]{
    \begin{figure}[htbp]                        %
        \centerline{\psfig{figure=#1/#2.ps,width=5in}}       %
        \small  %\vspace*{-2.7 in}                             %
        \begin{singlespace}                               %
        \caption{#3 \label{#2}}                           %
        \end{singlespace}                                 %
        \end{figure}                          %
    }                                         %

\newcommand{\etal}{{\sl et al.}}
\newcommand{\var}{\mbox{var}}
\newcommand{\corr}{\mbox{corr}}
\newcommand{\inte}{\int_{-\infty}^{+\infty}}
\newcommand{\btheta}{\theta} %\mbox{\boldmath $\theta$}}

\title{A selective overview of nonparametric methods in financial
econometrics}
\author{Jianqing Fan \\
Department of Operation Research and Financial Engineering\\
Princeton University\\
Princeton, NJ 08544
\thanks{Fan was partially supported by NSF grant DMS-0355179
    and a direct allocation RGC grant of the Chinese University of
    Hong Kong.  He acknowledges gratefully various discussions
    with Professors Yacine A\"it-Sahalia and Jia-an Yan and helpful
    comments of the editors and reviewers that lead to improve significantly
    the presentation of the paper.}
       }

\begin{document}
% typeset front matter
\maketitle

\begin{abstract}
This paper gives a brief overview on the nonparametric techniques
that are useful for financial econometric problems.  The problems
include estimation and inferences of instantaneous returns and
volatility functions of time-homogeneous and time-dependent
diffusion processes, and estimation of transition densities and
state price densities. We first briefly describe the problems and
then outline main techniques and main results. Some useful
probabilistic aspects of diffusion processes are also briefly
summarized to facilitate our presentation and applications.
\end{abstract}

\section{Introduction}

Technological invention and trade globalization have brought into
a new era of financial markets.  Over the last three decades,
enormous number of new financial products have been introduced to
meet customers' demands.  An important milestone is that in the
year 1973, the world's first options exchange opened in Chicago.
At the same year, Black and Scholes (1973) published their famous
paper on option pricing and Merton (1973) launched general
equilibrium model for security pricing, two landmarks for modern
asset pricing.  Since then, the derivative markets have
experienced extraordinary growth.  Professionals in finance now
routinely use sophisticated statistical techniques and modern
computation power in portfolio management, securities regulation,
proprietary trading, financial consulting and risk management.

Financial econometrics is an active field of integration of
finance, economics, probability, statistics, and applied
mathematics.  This is exemplified in the books by Campbell \etal
(1997), Gouri\'eroux and Jasiak (2001), and Cochrane (2001).
Financial activities generate many new problems, economics
provides useful theoretical foundation and guidance, and
quantitative methods such as statistics, probability and applied
mathematics are essential tools to solve the quantitative problems
in finance. To name a few, complex financial products pose new
challenges on their valuation and risk management. Sophisticated
stochastic models have been introduced to capture the salient
features of underlying economic variables and use for security
pricing. Statistical tools are used to identify parameters of
stochastic models, to simulate complex financial systems and to
test economic theories via empirical financial data.

An important area of financial econometrics is to study the
expected returns and volatilities of the price dynamics of stocks
and bonds. Returns and volatilities are directly related to asset
pricing, proprietary trading, security regulation and portfolio
management. To achieve these objectives, the stochastic dynamics
of underlying state variables should be correctly specified.  For
example, option pricing theory allows one to value stock or index
options and hedge against the risks of option writers, once a
model for the dynamics of underlying state variables is given.
See, for example, the books on mathematical finance by Bingham and
Kiesel (1998), Steele (2000), and Duffie (2001). Yet, many of
stochastic models in use are simple and convenient ones to
facilitate mathematical derivations and statistical inferences.
They are not derived from any economics theory and hence can not
be expected to fit all financial data. Thus, while the pricing
theory gives spectacularly beautiful formulas when the underlying
dynamics is correctly specified, it offers little guidance in
choosing or validating a model.  There is always a danger that
misspecification of a model leads to erroneous valuation and
hedging strategies. Hence, there are genuine needs for flexible
stochastic modeling. Nonparametric methods offer a unified and
elegant treatment for such a purpose.

Nonparametric approaches have recently been introduced to estimate
return, volatility, transition densities and state price densities
of stock prices and bond yields (interest rates). They are also
useful for examing the extent to which the dynamics of stock
prices and bond yields vary over time. They have immediate
applications to the valuation of bond price and stock options and
management of market risks. They can also be employed to test
economic theory such as the capital asset pricing model and
stochastic discount model (Campbell \etal 1997) and answer the
questions such as if the geometric Brownian motion fits certain
stock indices, whether the Cox-Ingsoll-Ross model fits yields of
bonds, and if interest rates dynamics evolve with time.
Furthermore, based on empirical data, one can also fit directly
the observed option prices with their associated characteristics
such as strike price, the time to maturity, risk-free interest
rate, dividend yield and see if the option prices are consistent
with the theoretical ones. Needless to say, nonparametric
techniques will play an increasingly important role in financial
econometrics, thanks to the availability of modern computing power
and the development of financial econometrics.

The paper is organized as follows.  We first introduce in section
2 some useful stochastic models for modeling stock prices and bond
yields and then briefly outline some probabilistic aspects of the
models.  In section 3, we review nonparametric techniques used for
estimating the drift and diffusion functions, based on either
discretely or continuously observed data.  In section 4, we
outline techniques for estimating state price densities and
transition densities. Their applications in asset pricing and
testing for parametric diffusion models are also introduced.
Section 5 makes some concluding remarks.

\section{Stochastic diffusion models}

Much of financial econometrics concerns with asset pricing,
portfolio choice and risk management.  Stochastic diffusion models
have been widely used for describing the dynamics of underlying
economic variables and asset prices.  They form the basis of many
spectacularly beautiful formulas for pricing contingent claims.
For an introduction to financial derivatives, see Hull (2003).

\subsection{One-factor diffusion models}

Let $S_{t\Delta}$ denote the stock price observed at time
$t\Delta$. The time unit can be hourly, daily, weekly, among
others. Presented in Figure 1(a) is the daily log-returns, defined
as
$$
  \log(S_{t\Delta}) - \log(S_{(t-1)\Delta}) \approx (S_{t\Delta} -
  S_{(t-1)\Delta})/S_{(t-1)\Delta}
$$
of the Standard and Poor 500 index, which is a value-weighted
index based on the prices of the 500 stocks that account for
approximately 70\% of the total U.S. equity (stock) market
capitalization. The styled features of the returns include that
the volatility tends to cluster and that the mean and variance of
the returns tend to be constant.  One simplified model to capture
the second feature is that
$$
    \log(S_{t\Delta}) - \log(S_{(t-1)\Delta}) \approx \mu_0 + \sigma_0
      \varepsilon_{t},
$$
where $\{\varepsilon_t\}$ is a sequence of independent normal
random variables.  This is basically a random walk hypothesis,
regarding the stock price movement as an independent random walk.
When the sampling time unit $\Delta$ gets small, the above random
walk can be regarded as a random sample from the continuous-time
process:
\begin{equation}
    d \log(S_t) = \mu_0 + \sigma_1 d W_t, \label{b1}
\end{equation}
where $\{W_t\}$ is a standard one-dimensional Brownian motion and
$\sigma_1 = \sigma_0 /\sqrt{\Delta}$.  The process (\ref{b1}) is
called geometric Brownian motion, as $S_t$ is an exponent of
Browian motion $W_t$.  It was used by Osborne (1959) to model
stock price dynamic and by Black and Scholes (1973) to derive
their celebrated option price formula.

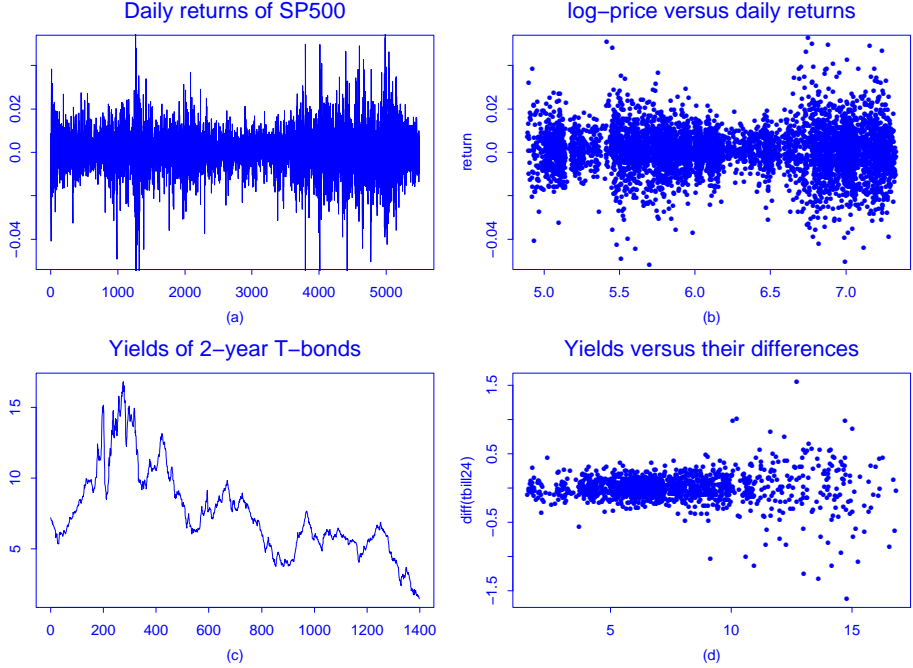
\fig{figs}{Fig1} {(a) Daily log-returns of the Standard \& Poor
500 index from October 21, 1980 to July 29, 2004. (b)  Scatter
plot of the returns against logarithm of the index (price level).
(c) Interest rates of two-year US Treasury notes from June 4, 1976
to March 7, 2007 sampled at weakly frequency.  (d)  Scatter plot
of the difference of yields versus the yields.}

Interest rates are fundamental to financial markets, consumer
spending, corporate earnings, asset pricing, inflation and
economy. The bond market is even bigger than the equity market.
Presented in Figure 1(c) is the interest rates $\{r_t\}$ of the
two-year US Treasury notes at weekly frequency.  As the interest
rates get higher, so do the volatilities.  To appreciate this,
Figure 1(d) plot the pairs $\{(r_{t-1}, r_{t} - r_{t-1})\}$. Its
dynamic is very different from the equity market.  The interest
rates should be non-negative.  They possess heteroscedasticity in
addition to the mean-revision property: As the interest rates rise
above the mean level $\alpha$, there is a negative drift that
pulls the rates down, while when the interest rates fall below
$\alpha$, there is a positive force that drives the rate up.  To
capture these two main features, Cox \etal (1985) derived the
following model for interest rate dynamic:
\begin{equation}
   dr_t = \kappa (\alpha - r_t) dt + \sigma
               r_t^{1/2} dW_t. \label{b2}
\end{equation}
For simplicity, we will refer to it as the CIR model. It is an
amelioration of the Vasicek (1977) model
\begin{equation}
  dr_t = \kappa (\alpha - r_t) dt + \sigma dW_t,
                             \label{b3}
\end{equation}
which ignores the heteroscedasticity and is also referred to as
the Ornstein-Uhlenback process. While this is an unrealistic model
for interest rates, the process is Gaussian with explicit
transition density.  It fact, the time series sampled from
(\ref{b3}) follows the autoregressive model of order 1:
\begin{equation}
   Y_t = (1 - \rho) \alpha + \rho Y_{t-1} + \varepsilon_t,
   \label{b4}
\end{equation}
where $Y_t = r_{t\Delta}$, $\varepsilon \sim N(0, \sigma^2
(1-\rho^2) /(2\kappa))$ and $\rho = \exp(-\kappa \Delta)$. Hence,
the process is well understood and usually serves as a testing
case for proposed statistical methods.

There are many stochastic models that have been introduced to
model the dynamics of stocks and bonds.  Let $X_t$ be an observed
economic variable at time $t$.  This can be the prices of a stock
or a stock index, or the yields of a bond. A simple and frequently
used stochastic model is
\begin{equation}
dX_t = \mu(X_t) dt + \sigma(X_t) dW_t. \label{b5}
\end{equation}
The function $\mu(\cdot)$ is frequently called a drift or
instantaneous return function and $\sigma(\cdot)$ is referred to
as a diffusion or volatility function, since
\begin{eqnarray*}
 \mu(X_t) =  \lim_{\Delta \to 0} \Delta^{-1} E(X_{t+\Delta} - X_t|X_t),
 \quad \and \quad
 \sigma^2(X_t) = \lim_{\Delta \to 0} \Delta^{-1} \var(X_{t+\Delta}|
 X_t).
\end{eqnarray*}

The time-homogeneous model (\ref{b5}) contains many famous
one-factor models in financial econometrics.  In an effort to
improve the flexibility of modeling interest dynamics, Chan \etal
(1992) extends the CIR model (\ref{b2}) to
\begin{eqnarray}
 dX_t = \kappa (\alpha - X_t) dt + \sigma X_t^\gamma dW_t.
                \label{b6}
\end{eqnarray}
A\"it-Sahalia (1996b) introduces nonlinear mean reversion: while
interest rates remain in the middle part of their domain, there is
little mean reversion and at the end of the domain, strong
nonlinear mean reversion emerges. He imposes the nonlinear drift
of form $(\alpha_0 X_t^{-1} + \alpha_1 + \alpha_2 X_t + \alpha_2
X_t^2)$. See also Ahn and Gao (1999), which models the interest
rates by $Y_t = X_t^{-1}$, in which the $X_t$ follows the CIR
model.

Economic conditions vary over time.  Thus, it is reasonable to
expect that the instantaneous return and volatility depend on both
time and price level for a given state variable, such as stock
prices and bond yields.  This leads to a further generalization of
model (\ref{b5}) to allow the coefficients to depend on time $t$:
\begin{equation}
    dX_t = \mu(t, X_t) dt + \sigma(t, X_t) dW_t. \label{b7}
\end{equation}
Since only a trajectory of the process is observed [see Figure
1(c)], there is no sufficient information to estimate the
bivariate functions in (\ref{b7}) without further restrictions. A
useful specification of model (\ref{b7}) is
\begin{equation}
    dX_t=\{\alpha_0(t)+\alpha_1(t) X_t \}\, dt +
        \beta_0(t) X_t^{\beta_1(t)} \, dW_t.  \label{b8}
\end{equation}
This is an extension of the CKLS model (\ref{b6}) by allowing the
coefficients to depend on time and was introduced and studied by
Fan \etal (2003). Model (\ref{b8}) includes many commonly-used
time-varying models for the yields of bonds, introduced by Ho and
Lee (1986), Hull and White (1990), Black, Derman and Toy (1990),
Black and Karasinski (1991), among others.
% They assume respectively the following forms:
%\begin{eqnarray*}
%\mbox{HL:}  && d X_t = \mu(t) \, dt + \sigma(t) \, dW_t, \\
%\mbox{HW:}  && d X_t = \{\alpha_0(t) + \alpha_1(t) X_t\} \, dt +
%    \sigma(t) {X_t}^i \, dW_t, \quad i=0 \ \mbox{or}\ 0.5, \\
%\mbox{BDT:} && dX_t=\{\alpha_1(t) X_t +\alpha_2(t) X_t \log (X_t)\}\,
%    dt+\beta _0(t) X_t \, dW_t, \\
%\mbox{BK:}  && dX_t=\{\alpha_1(t) X_t +\alpha_2(t) X_t \log (X_t)\}\,
%    dt+\beta _0(t) X_t \, dW_t, \ \mbox{with}\
%    \alpha_2(t)=\frac {d \log \{\beta_0(t)\}} {dt}.
%\end{eqnarray*}
The experience in Fan \etal (2003) and other studies of the
varying coefficient models (Chen and Tsay 1993, Hastie and
Tibshirani, 1993,  Cai \etal 2000) shows that coefficient
functions in (\ref{b8}) can not be estimated reliably due to the
collinearity effect in local estimation: localizing in the time
domain, the process $\{X_t\}$ is nearly constant. This leads Fan
\etal (2003) to introduce the semiparametric model:
\begin{equation}
dX_t = \{\alpha_0(t)+\alpha_1 X_t \}\, dt + \beta_0(t) X_t^{\beta}
\, dW_t.  \label{b9}
\end{equation}
to avoid the collinearity.

\subsection{Some probabilistic aspects}

A question arises naturally when there exists a solution to the
stochastic differential equation (SDE) (\ref{b7}).  Such a program
was first carried out by It\^o(1942, 1946).  For SDE (\ref{b7}),
there are two different meanings of solutions: strong solution and
weak solution. See sections 5.2  and 5.3 of Karatzas and Shreve
(1991).  Basically, for a given initial condition $\xi$, a strong
solution requires that $X_t$ is determined completely by the
information up to time $t$.
% namely, the information contained in the driving Brownian motion $\{W_s, 0
%\leq s \leq t\}$ and the initial condition $\xi$.
%The weak solutions relax the above requirement.
Under the Lipchitz and linear growth conditions on the drift and
diffusion functions,
%\[
% \| \mu(t, x) - \mu(t, y) \| + \| \sigma(t, x) - \sigma(t,y) \|
% \leq K \|x - y\|
%\]
%\[
%   \| \mu(t, x) \|^2 + \| \sigma(t, x) \|^2 \leq K ( 1 + \|x\|^2)
%\]
%for every $t, x$ and $y$,
for every $\xi$ that is independent of $\{W_s\}$, there exists a
strong solution of equation (\ref{b7}).  Such a solution is
unique.  See Theorem 2.9 of Karatzas and Shreve (1991).  %The above
%result holds for the process in $d$-dimensional Euclidean space.

%For time homogeneous model (\ref{b1}), Shorohod (1965) and Stroock and
%Varadhan (1969) proved that for every initial distribution $\xi$ such
%that $E\|\xi\|^{2m} < \infty$ for some $m > 1$, there exists a weak
%solution, provided that $\mu$ and $\sigma$ are bounded.
%See Theorem 5.4.22 (page 323) of Karatzas and Shreve (1991).

For one-dimensional time-homogeneous diffusion process (\ref{b5}),
weaker conditions can be obtained for the so-called weaker
solution.
%Let us recall the It\^o's formula: For process $X_t$ in
%(\ref{b12}), for a sufficiently regular function $f$ [see e.g.
%page 153 of Karatzas and Shreve (1991)],
%\[
%  df(X_t, t) = \left \{ \frac{\partial f(X_t, t)}{\partial t} +
%  \frac{1}{2} \frac{\partial^2 f(X_t, t)}{\partial x^2}
% \sigma^2(X_t, t) \right \} dt + \frac{\partial f(X_t, t)}{\partial x}
%  dX_t.
%\]
%By an application of It\^o's formula, for the process $X_t$ in
%(\ref{b1}), the process $Y_t = p(X_t)$ satisfies
%\[
%   dY_t = (\mu p' + \frac{1}{2} \sigma^2 p'') dt + \sigma p' d B_t.
%\]
%Hence, the drift term can be removed if one takes $p$ to satisfy
%\[
%   \mu p' + \frac{1}{2} \sigma^2 p'' = 0
%\]
%or
%\[
%   p'(x) = c \exp \left \{ -2 \int_0^x \mu(u)/\sigma^2(u) du \right \},
%\]
%for an arbitrary constant $c$.
%Note that if $\sigma$ is strictly
%positive and differentiable, the process (\ref{b1}) can also
%be transformed to have unit diffusion
%via $ Y_t = \int_\cdot^{X_t} \sigma(u)^{-1} du \equiv G(Y_t)$
%for some nonnegative lower limit, which satisfies
%\begin{equation}
%    dY_t = \mu_Y(Y_t) dt + dW_t,
%    \label{b14}
%\end{equation}
%where $\mu_Y(y) = (\mu/\sigma - \frac{1}{2} \sigma')(G^{-1}(y))$.
By an application of the It\^o's formula to an appropriate
transform of the process, one can make the drift to zero.   Thus,
we can consider without loss of generality that the drift in
(\ref{b5}) is zero.
%Let
%\[
%Z(\sigma) = \{x:  \sigma(x) = 0\} \mbox{ and  }
%I(\sigma) = \left \{ x:  \int_{x-}^{x+} \sigma(y)^{-2} dy = \infty \right \}.
%\]
%Note that if $\sigma(x)$ is continuous, then $I(\sigma) \subset
%N(\sigma)$. Engelbert and Schmidt (1984) prove that equation
%(\ref{b5}) with $\mu = 0$ has a weak solution for every initial
%distribution if and only if $I(\sigma) \subset N(\sigma)$. In that
%case, the uniqueness in law holds for every initial distribution
%if and only if $I(\sigma) = N(\sigma)$.
For such a model, Engelbert and Schmidt (1984) give a necessary
and sufficient condition of the existence of the solution. The
continuity of $\sigma$ suffices for the existence of the weak
solution. See Theorem 5.5.4 (page 333) of Karatzas and Shreve
(1991) and Theorem 23.1 of Kallenberg (2001).

We will again use several times the It\^{o} formula. For process
$X_t$ in (\ref{b7}), for a sufficiently regular function $f$ (
Karatzas and Shreve, 1991, p.153),
\begin{equation}
  df(X_t, t) = \left \{ \frac{\partial f(X_t, t)}{\partial t} +
  \frac{1}{2} \frac{\partial^2 f(X_t, t)}{\partial x^2}
 \sigma^2(X_t, t) \right \} dt + \frac{\partial f(X_t, t)}{\partial x}
  dX_t. \label{b10}
\end{equation}
The formula can be understood as the second order Taylor expansion
of $f(X_{t+\Delta}, t + \Delta) - f(X_t, t)$ by noticing that
$(X_{t+\Delta} - X_t)^2$ is approximately $\sigma^2(X_t, t)
\Delta$.

%\subsection{Markovian property and transition density}

The Markovian property plays an important role in statistical
inference. According to Theorem 5.4.20 of Karatzas and Shreve
(1991), the solution $X_t$ to equation (\ref{b5}) is Markovian,
provided that the coefficient functions $\mu$ and $\sigma$ are
bounded on compact subsets.  Let $p_\Delta(y|x)$ be the transition
density, the conditional density of $X_{t+\Delta} = y$ given $X_t
= x$. The transition density must satisfy the forward and backward
Kolmogorov equations
%\[
%  \frac{\partial}{\partial t} p_t(y|x) = \frac{1}{2}
%  \frac{\partial^2 }{\partial y^2} [\sigma^2(y)p_t(y|x)]
%  - \frac{\partial }{\partial y} [\mu(y) p_t(y|x)]
%\]
%and backward Kolmogorov equation
%\[
%\frac{\partial}{\partial t} p_t(y|x) = \frac{1}{2}
%  \frac{\partial^2 }{\partial x^2} [\sigma^2(x)p_t(y|x)]
%    + \frac{\partial }{\partial x} [\mu(x) p_t(y|x)].
%\]
(page 282, Karatzas and Shreve 1991).

Under the linear growth and Lipchitz's conditions, and additional
conditions on the boundary behavior of functions $\mu$ and
$\sigma$, the solution to equation (\ref{b1}) is positive and
ergodic.
%One sufficient condition is at the left boundary $0$,
%$\mu(x) \leq K
%x^{-1}$ for some $K > 1$ for $x$ in $(0, \varepsilon)$, at the right
%boundary $\infty$, $\mu(x) \leq - K x^\alpha$ for some $\alpha > 0$, and
%\[
%       \limsup_{x \to 0 \mbox{ \scriptsize or } x \to \infty} -
%       \{ \mu^2(x)/\sigma^2(x) + \mu'(x) \} /\sigma^2(x) < \infty.
%\]
%One sufficient condition is at the left boundary,
%$\int_{0+} s(x) dx = + \infty$ and at the right boundary, $\int^{+\infty}
%s(x) dx = + \infty$, and
%\[
%   \int_0^\infty \sigma^{-2}(x) / s(x) dx < \infty,
%\]
%where $s(x) = \exp(- 2 \int_.^x \mu(y) \sigma^{-2}(y) dy)$ is the
%derivative of the scale function of the diffusion process
%(\ref{b1}).  See \S5.7 of Rogers and Williams (1987).
%Let $f(x)$ be the invariant density.  %Multiplying both sides of
%the forward Kolmogorov equation by a factor of $f(x)$ and
%integrating it out with respect to $x$ and noticing that $\int
%p_t(y|x) f(x) dx = f(y)$, we have
%\[
%  \frac{1}{2} \frac{\partial ^2}{\partial y^2} [ \sigma^2(y) f(y)]
%  - \frac{\partial}{\partial y} [ \mu(y) f(y) ] = 0.
%\]
%Solving this differential equation yields
The invariant density is given by
\begin{equation}
  f(x) = 2 C_0 \sigma^{-2}(x) \exp( 2 \int_.^x  \mu(y) \sigma^{-2}(y) dy),
  \label{b11}
\end{equation}
where $C_0$ is a normalizing constant and the lower limit of the
integral does not matter.  If the initial distribution is taken
from the invariant density, then the process $\{X_t\}$ is
stationary with the marginal density $f$ and transition density
$p_\Delta$.

Let $H_t$ be the operator defined by
\begin{equation}
     (H_t g) (x) = E(g(X_t) | X_0 = x), \quad x \in R,
     \label{b12}
\end{equation}
where $f$ is a Borel measurable bounded function on $R$.
%The operator satisfies the Feller semigroup property:  for $s, t \geq 0$,
%$H_{s+t} = H_s H_t = H_t H_s$.
A stationary process $X_t$ is said to satisfy the condition $G_2(s,
\alpha)$ of Rosenblatt (1970) if there exists an $s$ such that
\[
  \|H_s\|_2^2 = \sup_{\{f: Ef(X)=0\}}  \frac{E (H_sf)^2(X)}{E f^2(X)}
  \leq \alpha^2 < 1,
\]
namely the operator is contractive.  As a consequence of the
semigroup ($H_{s+t} = H_s H_t$) and contraction properties, the
condition $G_2$ implies (Banon, 1977) that for any $t \in [0,
\infty)$, $ \|H_t \| \leq \alpha^{t/s-1}$. The latter implies, by
the Cauchy-Schwartz inequality, that
\begin{equation}
\rho(t) = \sup_{g_1, g_2} \corr(g_1(X_0), g_2(X_t)) \leq
\alpha^{t/s-1}, \label{b13}
\end{equation}
That is, the $\rho$-mixing coefficient decays exponentially fast.
Banon and Nguyen (1981) show further that for stationary Markov
process, $\rho(t) \to 0$ is equivalent to (\ref{b13}), namely,
$\rho$-mixing and geometric $\rho$-mixing are equivalent.

\subsection{Valuation of contingent claims}

An important application of SDE is the pricing of financial
derivatives, such as options and bonds.  It forms beautiful modern
asset pricing theory and provides useful guidance in practice.
Hull (2003) and Duffie (2001) offer very nice introduction to the
field.

The simplest financial derivative is the European call option. A
call option is the right to buy an asset at a certain price $K$
(strike price) before or at expiration time $T$. A put option
gives the right to sell an asset at a certain price K (strike
price) before or at expiration. European options allow option
holders to exercise only at maturity, while American options can
be exercised at any time before expiration.  Most stock options
are American, while options on stock indices are European.

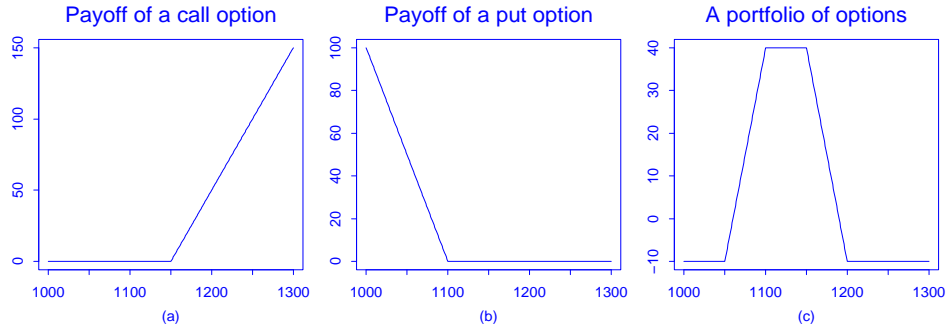
\fig{figs}{Fig2}{(a) Payoff of a call option.  (b)  Payoff of a
put option.  (c)  Payoff of a portfolio of 4 options with
different strike prices and different (long and short) positions}

The payoff for a European call option is $(X_T - K)_+$, where
$X_T$ is the price of the stock at expiration $T$.  When the stock
raises above the strike price $K$, one can excise the right and
makes a profit of $X_T - K$.  However, when the stock falls below
$K$, one renders his right and makes no profit.  Similarly, a
European put option has payoff $(K-X_T)_+$. See Figure 2. By
creating a portfolio with different maturity and different strike
prices, one can obtain all kind of payoff functions. As an
example, suppose that a portfolio of options consists of contracts
of SP500 index matured in 6 months:  one call-option with strike
price \$1,200, one put-option with strike price \$1,050, and \$40
cash, but short position (borrowing or $-1$ contract) on a call
option with strike price \$1,150 and a one put option with strike
\$1,100. Figure 2(c) shows the payoff function of such a portfolio
of options at the expiration $T$. Clearly, such an investor bets
the S\&P 500 index should be around \$1,125 in 6 months and limits
his risks on the investment. Thus, the European call and put
options are fundamental options as far as the payoff function at
time $T$ is concerned. There are many other exotic options such as
Asian options, look-back options and barrier options, which have
different payoff functions and the functions can be path
dependent. See Chapter 18 of Hull (2003).

Suppose that the asset price follows the SDE (\ref{b7}) and there
is a riskless investment alternative such as bond which earns
compounding rate of interest $r_t$.   Suppose that the underlying
asset pays no dividend.  Let $\beta_t$ be the value of the
riskless bond at time $t$. Then, with initial investment
$\beta_0$,
\[
   \beta_t = \beta_0 \exp( \int_0^t r_s ds),
\]
thanks to the compounding of interests. Suppose that the
probability measure $Q$ is equivalent to the original probability
measure $P$, namely $P(A)= 0$ if and only if $Q(A) = 0$.  The
measure $Q$ is called an equivalent martingale measure for
deflated price processes of given securities if these processes
are martingales with respect to $Q$.
%and if the Radon-Nikodym derivative $dQ/dP$ has finite variance.
An equivalent martingale measure is also
referred to as a ``risk-neutral'' measure if the deflater is the bond
price process.  See Chapter 6 of Duffie (2001).

When the markets are dynamically complete, the price of the
European option with payoff $\Psi(X_T)$ with initial price
$X_0=x_0$ is
\begin{equation}
   P_0 = \exp(-\int_0^T r_s ds) E^Q (\Psi(X_T)|X_0 =x_0), \label{b14}
\end{equation}
where $Q$ is the equivalent martingale measure for the deflated
price process $X_t / \beta_t$. Namely, it is the discounted value
of the expected payoff in the risk neutral world.  The formula is
derived by using the so-called relative pricing approach, which
values the price of the option from given prices of a portfolio,
consisting of the risk-free bond and the stock, with the identical
payoff as the option at the expiration.

As an illustrative example, suppose that the price of a stock
follows the geometric Brownian motion $dX_t = \mu X_t dt + \sigma
X_t dW_t$ and that the risk-free rate $r$ is constant. Then, the
deflated  price process $Y_t = \exp(-rt) X_t$ follows the SDE:
$$
   dY_t = (\mu - r) Y_t dt + \sigma Y_t d W_t.
$$
The deflated price process is not a martingale as the drift is not
zero. The risk-neutral measure is the one that makes the drift
zero.  To achieve this, we appeal to the Girsanov theorem, which
changes the drift of a diffusion process without alternating the
diffusion, via a change of probability measure.
%With the probability measure $Q$ defined as
%$$
%   \frac{dQ}{dP} = \exp (- \int_0^T \theta_s dW_s - 0.5 \int_0^T
%    \theta_s ds ),
%$$
%where $\theta_s = (\mu - r)/\sigma$,
Under the ``risk-neutral'' probability measure $Q$, the process
$Y_t$ satisfies $dY_t = \sigma Y_t dW_t$, a martingale. Hence, the
price process $X_t = \exp(rt) Y_t$  under $Q$ follows
\begin{equation}
    dX_t = r X_t dt + \sigma X_t dW_t. \label{b15}
\end{equation}
Using exactly the same derivation, one can easily generalize the
result to the price process (\ref{b5}).  Under the risk-neutral
measure, the price process (\ref{b5}) follows
\[
    dX_t = r X_t dt + \sigma (X_t) dW_t.
\]
The intuitive explanation of this is clear: all stocks under the
``risk-neutral'' world is expected to earn the same rate as the
risk-free bond.

For the Geometric Brownian motion, by an application of the It\^o
formula (\ref{b10}) to (\ref{b15}), we have under the
``risk-neutral'' measure
\begin{equation}
    \log X_t - \log X_0 = (r - \sigma^2/2) t + \sigma^2 W_t. \label{b16}
\end{equation}
Note that given the initial price $X_0$, the price follows a
log-normal distribution.  Evaluation of the expectation of
(\ref{b14}) for the European call option with payoff $\Psi(X_T) =
(X_T - K)_+$, one obtains the Black-Scholes (1973) option pricing
formula:
\begin{equation}
P_0 = x_0 \Phi(d_1) - K\exp(-rT)\Phi(d_{2}), \label{b17}
\end{equation}
where $d_{1} = \{\log(x_{0}/K) + (r + \sigma^{2}/2)T\}\{\sigma
\sqrt{T}\}^{-1}$ and $d_{2} = d_{1} - \sigma\sqrt{T}$.

\subsection{Simulation of stochastic models}

Simulation methods provide useful tools for valuation of financial
derivatives and other financial instruments, when analytic formula
(\ref{b14}) is hard to obtain.  They also provide useful tools for
assessing performance of statistical methods and statistical
inferences.

The simplest method is perhaps the Euler scheme.  The SDE
(\ref{b7}) is approximated as
\begin{equation}
    X_{t+\Delta} = X_t + \mu(t, X_t) \Delta + \sigma(t, X_t)
    \Delta^{1/2} \varepsilon_t,  \label{b18}
\end{equation}
where $\{\varepsilon_t\}$ is a sequence of independent random
variables with the standard normal distribution.  The time unit is
usually a year. Thus, the monthly, weekly and daily data
correspond, respectively, to $\Delta = 1/12, 1/52$ and $1/252$
(there are approximately 252 trading days per year).  Given an
initial value, one can recursively apply (\ref{b18}) to obtain a
sequence of simulated data $\{X_{j\Delta}, j=1,2, \cdots\}$.  The
approximation error can be reduced if one uses a smaller step size
$\Delta / M$ for a given integer $M$ to obtain first a more
detailed sequence $\{X_{j\Delta/M}, j=1, 2, \cdots \}$ and then to
take the subsequence $\{X_{j\Delta}, j= 1, 2, \cdots\}$. For
example, to simulate daily prices of a stock, one can simulate
hourly data first and than takes the daily closing prices. Since
the step size $\Delta/M$ is smaller, the approximation (\ref{b18})
is more accurate. However, the computational cost is about a
factor of $M$ higher.

The Euler scheme has convergence rate $\Delta^{1/2}$, which is
called strong order $0.5$ approximation by Kloeden \etal (1996).
The higher order approximations can be obtained by the
It\^o-Taylor expansion (see Schurz, 2000, page 242).  In
particular, a strong order-one approximation is given by
\begin{equation}
    X_{t+\Delta} = X_t + \mu(t, X_t) \Delta + \sigma(t, X_t)
    \Delta^{1/2} \varepsilon_t + \frac{1}{2} \sigma (t, X_t)
    \sigma_x'(t, X_t) \Delta \{ \varepsilon_t^2 -1\},  \label{b19}
\end{equation}
where $\sigma_x'(t,x)$ is the partial derivative function with
respect to $x$.  This method can be combined with a smaller step
size method in the last paragraph. For the time-homogeneous model
(\ref{b1}), an alternative form, without evaluating the derivative
function, is
%\[
%    X_{t+\Delta} = X_t + \mu(X_t) \Delta + \frac{1}{2} \{
%    \sigma(X_t) + \sigma(X_t+\mu(X_t) \Delta + \sigma(X_t)
%    \Delta^{1/2} \varepsilon_t) \}
%    \Delta^{1/2} \varepsilon_t.
%\]
%See
given in (3.14) of Kloeden \etal (1996).
%  More numerical methods can be found in Kloeden and Platen (1995).

The exact simulation method is available if one can simulate the
data from the transition density.  Given the current value
$X_t=x_0$, one draws $X_{t+\Delta}$ from the transition density
$p_\Delta(\cdot|x_0)$. The initial condition can either be fixed
at a given value or be generated from the invariant density
(\ref{b11}).  In the latter case, the generated sequence is
stationary.

There are only a few processes where exact simulation is possible.
For GBM, one can generate the sequence from the explicit solution
(\ref{b16}), where the Brownian motion can be simulated from
independent Gaussian increments.  The conditional density of
Vascicek's model (\ref{b3}) is Gaussian with mean $\alpha + (x_0 -
\alpha) \rho$ and variance $\sigma^2_{\Delta} = \sigma^2 (1 -
\rho^2) /(2\kappa)$, as indicated by (\ref{b4}).
%Indeed, the solution to the SDE (\ref{b6}) is
%\[
%  X_t = \alpha + (X_0 - \alpha) \exp( -\kappa t) + \sigma \exp(-
%  \kappa t) \int_0^t \exp(\kappa s) dW_s.
%\]
%By (\ref{b11}), one can easily see that
%its invariant density is Gaussian with mean $\alpha$ and variance
%$\sigma^2/(2\kappa)$.
Generate $X_0$ from the invariant density $N(\alpha,
\sigma^2/(2\kappa))$. With $X_0$, generate $X_\Delta$ from the
normal distribution with mean $\alpha + (X_0 - \alpha)
\exp(-\kappa \Delta)$ and variance $\sigma_\Delta^2$. With
$X_\Delta$, we generate $X_{2\Delta}$ from $\alpha + (X_\Delta -
\alpha) \exp(-\kappa \Delta)$ and variance $\sigma_\Delta^2$.
Repeat this process until we obtain the desired length of the
process.

For the CIR model (\ref{b2}), provided that $q = 2 \kappa \alpha
/\sigma^2 - 1 \geq 0$ (a sufficient condition for $X_t \geq 0$),
%the transition density 
%is given by $p_\Delta(y|x) = c \exp(-u -v)
%$I_q(\cdot???) (v/u)^{q/2}$ with $c = 2 \kappa /\{\sigma^2 (1 -
%$\exp(-\kappa \Delta))\}$, $u = cx_0 \exp(k \Delta)$, $v = cx$ and
%$$
%I_q(x) = \left ( \frac{x}{2} \right )^q \sum_{i=0}^\infty
%\frac{1}{i! \Gamma(q+i+1)} \left ( \frac{x}{2} \right )^{2i}
%$$
%is the modified Bessel function of the first kind of order $q$.  
the transition density can be expressed in terms of the modified Bessel
function of the first kind.  This distribution is often referred to
as the noncentral
$\chi^2$ distribution.  That is, given $X_t = x_0$, $2cX_{t+\Delta}$
has a noncentral $\chi^2$ distribution with degrees of freedom $2q+2$
and noncentrality parameter $2u$.
The invariant density is the Gamma distribution with shape parameter
$q+1$ and the scale parameter $\sigma^2/(2 \kappa)$.

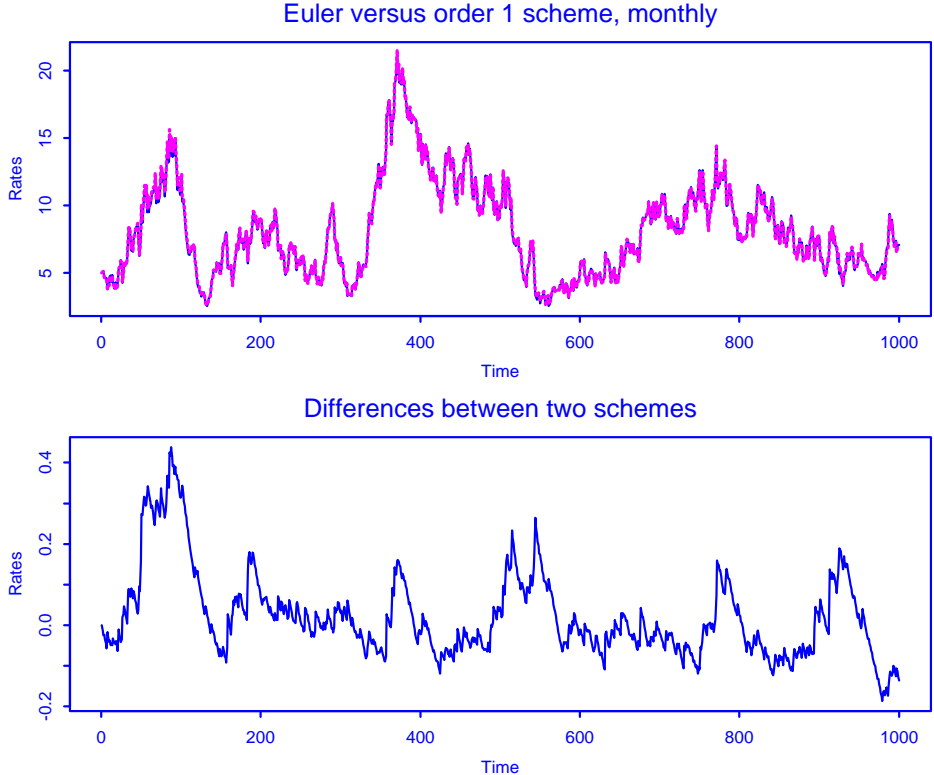
\fig{figs}{Fig3} {Simulated trajectories (multiplied by 100) using
the Euler approximation and the strong order-one approximation for
a CIR model. Top panel:  solid-curve corresponds to the Euler
approximation and the dashed curve is based on the order-one
approximation.  Botton panel:  The difference between the
order-one scheme and the Euler scheme.}

As an illustration, we consider the CIR model (\ref{b7}) with
parameters $\kappa = 0.21459$, $\alpha = 0.08571$, $\sigma =
0.07830$ and $\Delta = 1/12$.  The model parameters are taken from
Chapman and Pearson (2000). We simulated 1000 monthly data using
both the Euler scheme (\ref{b18}) and strong order-one
approximation (\ref{b19}) with the same random shocks.  Figure 3
depicts one of their trajectories.  The difference is negligible.
This is in line with the observations made by Stanton (1997) that
as long as data are sampled monthly or more frequently, the errors
introduced by using the Euler approximation  is very small for
stochastic dynamics that are similar to the CIR model.

\section{Estimation of return and volatility functions}

There are a large literature on the estimation of the return and
volatility functions.  Early references include Pham (1981) and Prakasa
Rao (1985).  Some studies are based on continuously observed
data, while others are based on discretely observed data.  For the
latter, some regard $\Delta$ tending to zero, while others regard
$\Delta$ fixed.   We briefly introduce some of the ideas.

\subsection{Methods of estimation}

We first outline several methods of estimation for parametric
models. The idea can be expanded into nonparametric models.
Suppose that we have a sample $\{X_{i\Delta}, i =0, \cdots, n \}$
from model (\ref{b5}). Then, the likelihood function, under the
stationary condition, is
\begin{equation}
   \log f(X_0) + \sum_{i=1}^n \log  p_\Delta (X_{i\Delta}|X_{(i-1)\Delta}).
   \label{c1}
\end{equation}
If the functions $\mu$ and $\sigma$ are parameterized and the
explicit form of the transition density is available, one can
apply the maximum likelihood method.   However, the explicit form
of the transition density is not available for many simple models
such as the CLKS model (\ref{b6}).  Even for the CIR model
(\ref{b2}), its maximum likelihood estimator is very difficult to
find.

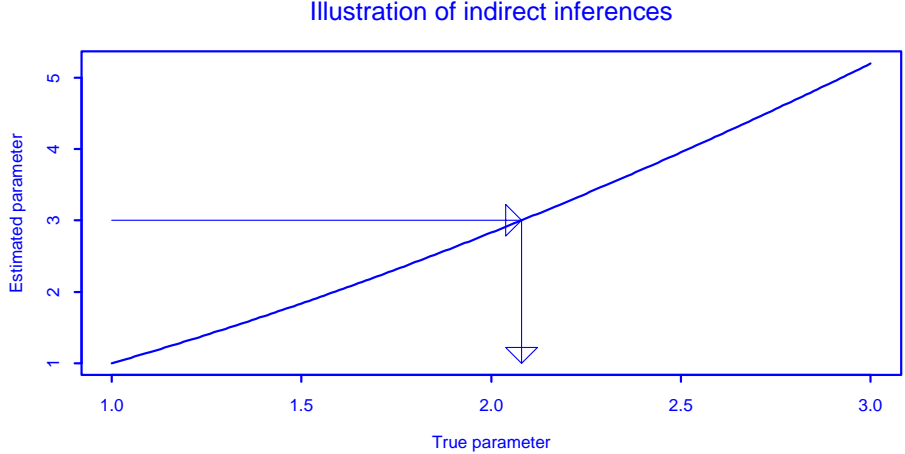
\fig{figs}{Fig4}{Illustration the idea of the indirect inference.
For each given true $\btheta$, one obtains an estimate using the
Euler approximation.  This gives a calibration curve as shown. Now
for a given estimate $\hat{\btheta}_0 = 3$ based on the Euler
approximation, one find the calibrated estimate
$\hat{\btheta}_1^{-1}(3) = 2.080$.}

One simple technique is to rely on the Euler approximation scheme
(\ref{b18}).  Then proceed as if the data come from the Gaussian
location and scale model.  This method works well when $\Delta$ is
small, but can create some biases when $\Delta$ is large. However,
the bias can be reduced by the following calibration idea, called
the indirect inference by Gouri\'eroux \etal (1993).  The idea
works as follows. Suppose that the functions $\mu$ and $\sigma$
have been parameterized with unknown parameters $\btheta$.  Use
the Euler approximation (\ref{b18}) and the maximum likelihood
method to obtain an estimate $\hat{\btheta}_0$. For each given
parameter $\btheta$ around $\hat{\btheta}_0$ , simulate data from
(\ref{b5}) and apply the crude method to obtain an estimate
$\hat{\btheta}_1(\btheta)$, which depends on $\btheta$. Since we
simulated the data with the true parameter $\btheta$, the function
$\hat{\btheta}_1(\btheta)$ tells us how to calibrate the estimate.
See Figure 4. Calibrate the estimate via $\hat{\btheta}_1^{-1}
(\hat{\btheta}_0)$, which improves the bias of estimate.  One
drawback of this method is that it is intensive in computation and
the calibration can not easily be done when the dimensionality of
parameters $\theta$ is high.

Another method for bias reduction is to approximate the transition
density in (\ref{c1}) by a higher order approximation, and to then
maximize the approximated likelihood function. Such a scheme has
been introduced by A\"it-Sahalia (1999, 2002) who derives the
expansion of the transition density around a normal density
function using the Hermit polynomial.  The intuition behind such
an expansion is that the diffusion process $X_{t+ \Delta} - X_t$
in (\ref{b5}) can be regarded as sum of many independent
increments with a smaller step size and hence the Edgeworth
expansion can be obtained for the distribution of $X_{t+\Delta} -
X_t$ given $X_t$.

An ``exact'' approach is to use the method of moment.  If the
process $X_t$ is stationary as in the interest-rate models, the
moment conditions can easily be derived by observing
$$
E \{ \lim_{\Delta \to 0}  \Delta^{-1} E [ g(X_{t+\Delta}) - g(X_t)
| X_t ]\} =  \lim_{\Delta \to 0}  \Delta^{-1}  E [ g(X_{t+\Delta})
- g(X_t) ] = 0
$$
for any function $g$ satisfying the regularity condition such that
the limit and the expectation is exchangeable. The right-hand side
is the expectation of $d g(X_t)$.  By It\^o's formula (\ref{b10}),
the above equation deduces to
\begin{equation}
        E [ g'(X_t) \mu(X_t) + g''(X_t) \sigma^2(X_t)/2] = 0.
        \label{c2}
\end{equation}
For example, if $g(x) = \exp(-ax)$ for some given $a > 0$, then
$$
    E \exp(- a X_t)\{ \mu(X_t) - a \sigma^2(X_t)/2 \} = 0.
$$
This can produce arbitrary number of equations by choosing
different $a$'s.  If the functions $\mu$ and $\sigma$ are
parametrized, the number of moment conditions can be more than the
number of equations. One way to efficiently use this is the
generalized method of moment introduced by Hansen (1982),
minimizing a quadratic form of the discrepancies between the
empirical and the theoretical moments, a generalization of the
classical method of moment which solves the moment equations. The
weighting matrix in the quadratic form can be chosen to optimize
the performance of the resulting estimator. To improve the
efficiency of the estimate, a large system of moments is needed.
Thus, the generalized method of moments needs a large system of
nonlinear equations, which can be expensive in computation.
Further, the moment equations (\ref{c2}) use only the marginal
information of the process.  Hence, it is not efficient. For
example, in the CKLS model (\ref{b6}), $\sigma$ and $\kappa$ are
estimable via (\ref{c2}) only through $\sigma^2/\kappa$.

\subsection{Time-homogeneous model}

The Euler approximation can easily be used to estimate the drift
and diffusion nonparametrically. Let $Y_{i \Delta} = \Delta^{-1}
(X_{(i+1) \Delta} - X_{i\Delta})$ and $Z_{i\Delta} = \Delta^{-1}
(X_{(i+1) \Delta} - X_{i\Delta})^2$. Then,
\[
   E(Y_{i \Delta}|X_{i\Delta}) = \mu(X_{i\Delta}) + O(\Delta),
   \mbox{ and }
   E(Z_{i \Delta}|X_{i\Delta}) = \sigma^2(X_{i\Delta}) + O(\Delta).
\]
Thus, $\mu(\cdot)$ and $\sigma^2(\cdot)$ can be approximately
regarded as the regression functions of $Y_{i\Delta}$  and
$Z_{i\Delta}$ on $X_{i\Delta}$, respectively.  Stanton (1997)
applies kernel regression (Wand and Jones, 1995; Simonoff, 1996)
to estimate the return and volatility functions.   Let $K(\cdot)$
be a kernel function and $h$ be a bandwidth.  Stanton's estimators
are given by
\[
   \hat{\mu}(x) = \frac{\sum_{i=0}^{n-1} Y_{i \Delta} K_h(X_{i\Delta} - x) }
         { \sum_{i=0}^{n-1}  K_h(X_{i\Delta} - x)},
    \quad \mbox{and} \quad
   \hat{\sigma}^2(x) = \frac{\sum_{i=0}^{n-1}
       Z_{i \Delta} K_h(X_{i\Delta} - x) }
         { \sum_{i=0}^{n-1}  K_h(X_{i\Delta} - x)},
\]
where $K_h(u) = h^{-1} K(u/h)$ is a rescaled kernel. The
consistency and asymptotic normality of the estimator are studied
in Bandi and Phillips (1998). Independently, Fan and Yao (1998)
apply the local linear technique (\S6.3, Fan and Yao 2003) to
estimate the return and volatility functions, under a slightly
different setup.   The local linear estimator (Fan, 1992) is given
by
\begin{equation}
\hat{\mu}(x) = \sum_{i=0}^{n-1} K_n (X_{i\Delta}-x, x) Y_{i\Delta}
, \label{c3}
\end{equation}
where
\begin{equation}
    K_n(u, x) = K_h(u) \frac{S_{n,2}(x) - u S_{n,1}(x)}{S_{n,2}(x)
    S_{n,0}(x) - S_{n,1}(x)^2} \label{c4}
\end{equation}
with $S_{n,j}(x) = \sum_{i=0}^{n-1} K_h(X_{i\Delta} -
x)(X_{i\Delta} - x)^j$ is the equivalent kernel induced by the
local linear fit.   In contrast with the kernel method, the local
linear weights depend on both $X_i$ and $x$.  In particular, they
satisfy that
$$
  \sum_{i=1}^{n-1} K_n (X_{i\Delta}-x, x) = 1, \quad  \mbox{and} \quad
  \sum_{i=1}^{n-1} K_n (X_{i\Delta}-x, x) (X_{i\Delta}-x) = 0.
$$
These are the key properties for the bias reduction of the local
linear method as demonstrated in Fan (1992). Further, Fan and Yao
(1998) use the squared residuals
\[
  \Delta^{-1} (X_{(i+1) \Delta} - X_{i\Delta} - \hat{\mu}(X_{i\Delta})\Delta)^2
\]
rather than $Z_{i\Delta}$ to estimate the volatility function.
This would reduce further the approximation errors in the
volatility estimation.  They show further that the conditional
variance function can be estimated as well as if the conditional
mean function is known in advance.

\begin{table}[bhtp]
\centering\small \caption{\textsl{Variance inflation factors by
using higher order differences}}
\begin{center}
%\doublerulesep 0.5pt
\begin{tabular}{c r r r r r} \hline \hline \label{table1}
       & \multicolumn{5}{c}{\textsl{Order $k$}} \\  \cline{2-6}
       & \multicolumn{1}{c}{$1$} & \multicolumn{1}{c}{$2$}
       & \multicolumn{1}{c}{$3$} & \multicolumn{1}{c}{$4$}
       & \multicolumn{1}{c}{$5$}       \\ \hline
$V_1(k)$ & 1.00 & 2.50 & 4.83 & 9.25 & 18.95  \\
$V_2(k)$ & 1.00 & 3.00 & 8.00 &21.66 & 61.50  \\ \hline
\end{tabular}
\end{center}
\end{table}

Stanton (1997) derives higher order approximation scheme up to
order 3 in an effort to reduce biases. He suggested that higher
order approximations must outperform lower order approximations.
To verify such a claim, Fan and Zhang (2003) derive the following
order $k$ approximation scheme:
\begin{equation}
    E(Y_{i \Delta}^* |X_{i\Delta}) = \mu(X_{i\Delta}) + O(\Delta^k),
       \mbox{ and }
    E(Z_{i \Delta}^* |X_{i\Delta}) = \sigma^2(X_{i\Delta}) +
      O(\Delta^k),   \label{c5}
\end{equation}
where
\[
    Y_{i \Delta}^*  = \Delta^{-1} \sum _{j=1}^k a_{k,j}
        \{X_{(i+j)\Delta}-X_{i\Delta}\}
        \mbox{ and }
    Z_{i \Delta}^*  = \Delta^{-1} \sum _{j=1}^k a_{k,j}
        \{X_{(i+j)\Delta}-X_{i\Delta}\}^2
\]
and the coefficients $a_{k,j}=(-1)^{j+1} {k \choose j} \big/ j$
are chosen to make the approximation error in (\ref{c5}) of order
$\Delta^k$.  For example, the second approximation is
$$
 1.5 (X_{t+\Delta} - X_t) - 0.5 (X_{t+ 2 \Delta} - X_{t+\Delta}).
$$
By using the independent increment of the Brownian motion, its
variance is $1.5^2+0.5^2=2.5$ times as large as that of the first
order difference.  Indeed, Fan and Zhang (2003) show that while
higher order approximations give better approximation errors, we
have to pay a huge premium for variance inflation:
\begin{eqnarray*}
\var(Y_{i\Delta}^*|X_{i\Delta}) & = & \sigma^2(X_{i\Delta}) V_1(k)
    \Delta^{-1} \{1+O(\Delta)\}, \\
\var(Z_{i\Delta}^*|X_{i\Delta}) & = & 2\sigma^4(X_{i\Delta})
V_2(k) \{1+O(\Delta)\},
\end{eqnarray*}
where the variance inflation factors $V_1(k)$ and $V_2(k)$ are
explicitly given by Fan and Zhang (2003).  Table 1 depicts some of
the numerical results for the variance inflation factor.

The above theoretical results have also been verified via
empirical simulations in Fan and Zhang (2003). The problem is no
monopoly for nonparametric fitting --- it is shared by the
parametric methods.  Therefore, the methods based on higher order
differences should seldomly be used unless the sampling interval
is very wide (e.g. quarterly data). It remains open whether it is
possible to estimate nonparametrically the return and the
volatility functions, without seriously inflating the variance,
with other higher approximation schemes.

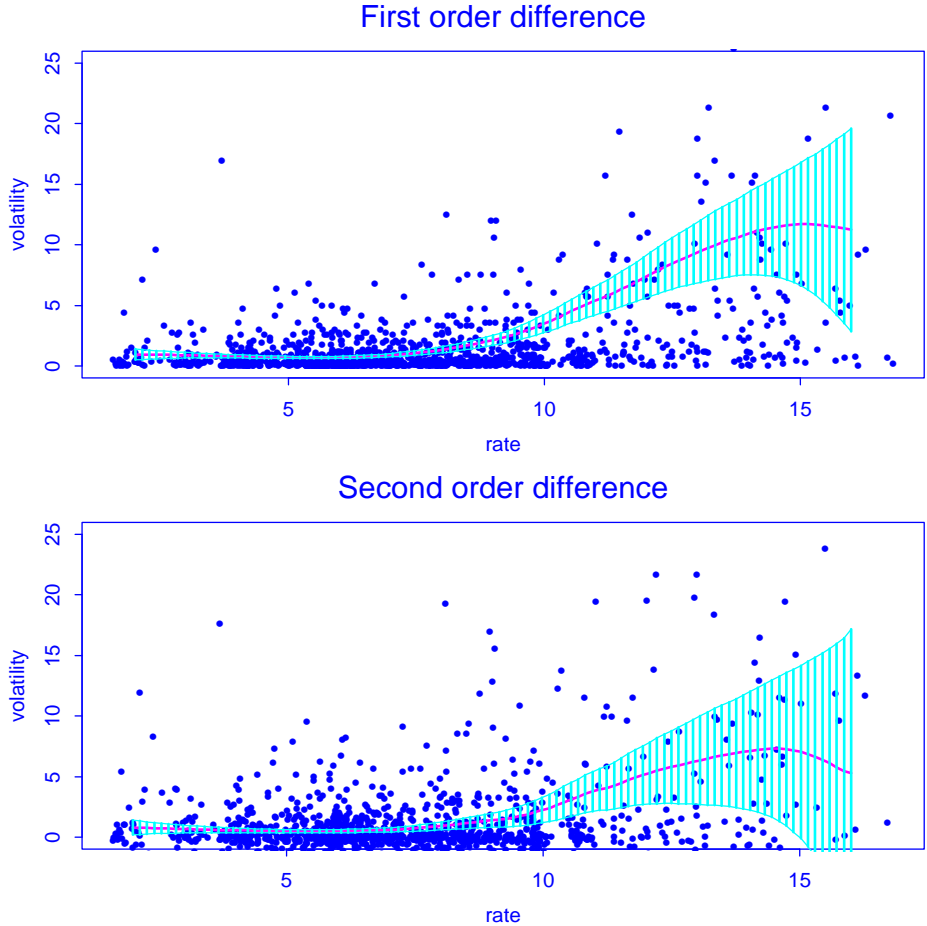
\fig{figs}{Fig5} {Nonparametric estimates of volatility based on
orders 1 and 2 differences.  The bars represent two standard
deviations above and below the estimated volatility.  Top panel:
order 1 fit.  Botton panel: order 2 fit.}

As an illustration, we take the yields of the two-year Treasury
notes depicted in Figure 1. Figure 5 presents nonparametrically
estimated volatility function based orders $k=1$ and $k=2$
approximations. The local linear fit is employed with the
Epanechnikov kernel and bandwidth $h = 0.35$. It is evident that
the order 2 approximation has higher variance than the order 1
approximation.  In fact, the magnitude of variance inflation is in
line with the theoretical result: the increase of the standard
deviation from order 1 to order 2 approximation is $\sqrt{3}$.

Stanton (1997) applies his kernel estimator to a Treasury's bill
data set and observes nonlinear return function in his
nonparametric estimate, particularly in the region where the
interest rate is high (over 14\%, say). This leads him to
postulate the hypothesis that the return functions of short-term
rates are nonlinear.  Chapman and Pearson (2000) study the finite
sample properties of Stanton's estimator.  By applying his
procedure to the CIR model, they find that the Stanton's procedure
produces spurious nonlinearity, due to the boundary effect and the
mean reversion.

Can we employ a formal statistic test to the Stanton's hypothesis?
The null hypothesis can simply be formulated:  the drift is of a
linear form as in model (\ref{b6}).  What is the alternative
hypothesis? For such a kind of problem, our alternative model is
usually vague. Hence, it is natural to assume that the drift is a
nonlinear smooth function.  This becomes a testing problem with a
parametric null hypothesis versus a nonparametric alternative
hypothesis. There is a large literature on this.  The basic idea
is to compute a discrepancy measure between the parametric
estimates and nonparametric estimates and to reject the parametric
hypothesis when the discrepancy is large.  See for example the
book by Hart (1997).  In an effort to derive a generally
applicable principle, Fan \etal (2001) propose the generalized
likelihood ratio (GLR) tests for parametric versus nonparametric
or nonparametric versus nonparametric hypotheses. The basic idea
is to replace the maximum likelihood under nonparametric
hypotheses (usually does not exist) by the likelihood under good
nonparametric estimates. The method has been successfully employed
by Fan and Zhang (2003) for checking whether the return and
volatility functions possess certain parametric forms.

Various discretization schemes and estimation methods have been
proposed for the case with high frequency data over a long time
horizon.  More precisely, the studies are under the assumptions
that $\Delta_n \to 0$ and $n \Delta_n \to \infty$.  See for
example, Dacunha-Castelle and Florens (1986), Yoshida (1992),
Kessler (1997), Arfi (1998), Gobet (2003), Cai and Hong (2003) and
references therein. Arapis and Gao (2003) investigate the mean
integrated square errors of several methods for estimating the
drift and diffusion and compare their performance. A\"it-Sahalia
and Mykland (2003, 2004) study the effects of random and discrete
sampling when estimating continuous-time diffusions. Bandi and
Nguyen (1999) investigate small sample behaviors of nonparametric
diffusion estimators. Thorough study of nonparametric estimation
of conditional variance functions can be found in M\"uller and
Stadtm\"uller (1987), Hall and Carroll (1989), Ruppert \etal
(1997) and H\"ardle and Tsybakov (1997). In particular, \S8.7 of
Fan and Yao (2003) give various methods for estimating the
conditional variance function.

%There is a large statistical literature on estimating the mean and
%variance functions.  See the books by Wand and Jones (1995), Fan
%and Gijbels (1996), Simonoff (1996) and Fan and Yao (2003).  In
%particular, \S8.7 of Fan and Yao (2003) give various methods for
%estimating the conditional variance function.  For the estimation
%of constant variance with nonparametric regression function, Rice
%(1984), Gasser \etal (1986) and Hall \etal (1990) proposed various
%root-n consistent estimators.  For the case where the conditional
%variance function is not constant, see M\"uller and Stadtm\"uller
%(1987, 1993), Hall and Carroll (1989), Ruppert \etal (1997) and
%H\"ardle and Tsybakov (1997).

\subsection{Fixed sampling interval}

For practical analysis of financial data, it is hard to determine
whether the sampling interval tends to zero.  The key
determination is whether the approximation errors for small
``$\Delta$'' are negligible.  It is ideal when a method is
applicable whether or not ``$\Delta$'' is small.  This kind of
method is possible, as demonstrated below.

The simplest problem to illustrate the idea is the kernel density
estimation of the invariant density of the stationary process
$\{X_t\}$.  For the given sample $\{X_{t\Delta} \}$, the kernel
density estimate for the invariant density is
\begin{equation}
    \hat{f}(x) = n^{-1}  \sum_{i=1}^n K_h (X_{i\Delta} - x),
    \label{c6}
\end{equation}
based on the discrete data $\{X_{i\Delta}, i=1, \cdots, n\}$. This
method is valid for all $\Delta$.  It gives a consistent estimate
of $f$ as long as the time horizon is long: $n\Delta \to \infty$.
We will refer to this kind of nonparametric methods as the
state-domain smoothing, as the procedure localizes in the state
variable $X_t$.  Various properties, including consistency and
asymptotic normality, of the kernel estimator (\ref{c6}), are
studied by Bandi (1998), Bandi and Phillips (1998).  Bandi (1998)
also used the estimator (\ref{c6}), which is the same as the local
time of the process spending at a point $x$ except a scaling
constant, as a descriptive tool for potentially nonstationary
diffusion processes.

Why can the state-domain smoothing methods be employed as if the
data were independent?  This is due to the fact that localizing in
state-domain weakens the correlation structure and that
nonparametric estimates use essentially only local data.  Hence
many results on nonparametric estimators for independent data
continue to hold for dependent data, as long as their mixing
coefficients decay sufficiently fast.  As mentioned at the end of
\S 2.2, the geometric mixing and mixing are equivalent for
time-homogeneous diffusion process.  Hence, the mixing
coefficients decay sufficiently fast for theoretical
investigation.

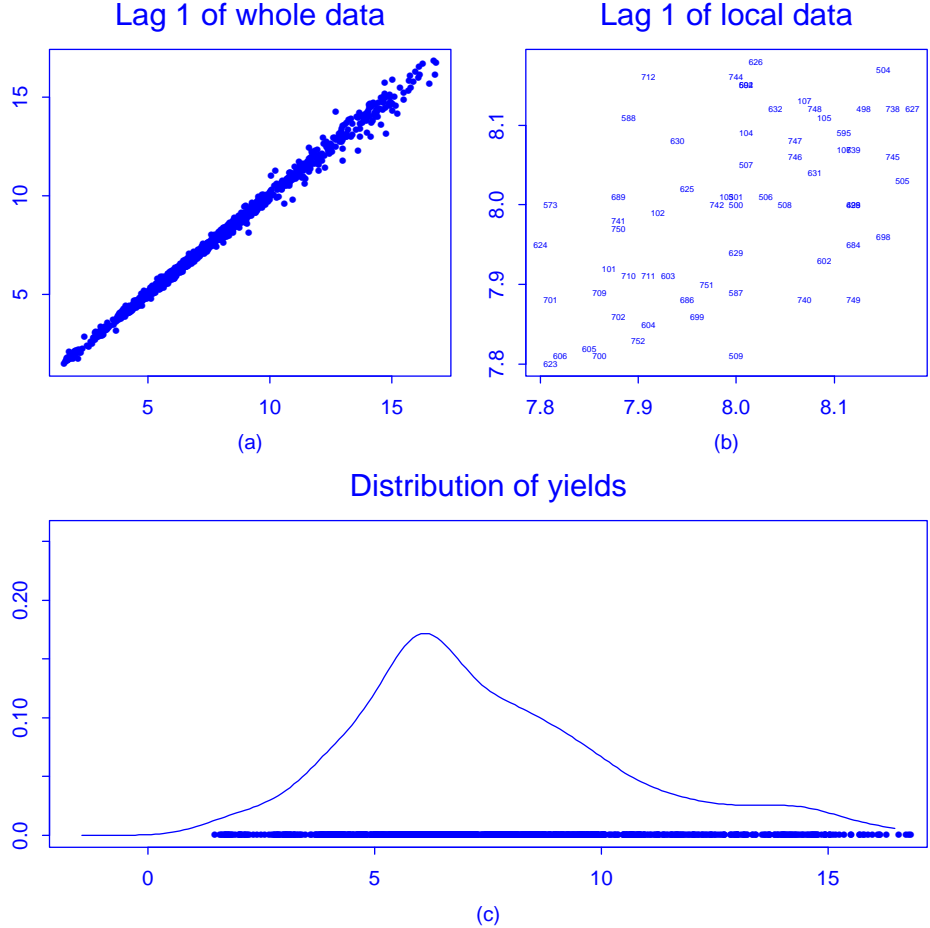
\fig{figs}{Fig6} {(a) Lag 1 scatterplot of the two-year Treasury
note data.
 (b)  Lag 1 scatterplot of those data falling in the neighborhood
      $8\% \pm 0.2\%$ --- the points are represented by the time of the observed
      data.  The number in the scatterplot shows the indices
      of the data falling in the neighborhood.  (c)  Kernel density
      estimate of the invariant density.
      }

The localizing and whitening can be understood graphically in
Figure 6. Figure 6(a) shows that there is very strong serial
correlation of the yields of the two-year treasury notes. However,
this correlation is significantly weakened for the local data in
the neighborhood of $8\% \pm 0.2\%$.  In fact, as detailed in
Figure 6(b), the indices of the data that fall in the local window
are quite far apart. This in turn implies the week dependence for
the data in the local window, i.e. ``whitening by windowing''. See
\S5.4 of Fan and Yao (2003) and Hart (1996) for further details.
The effect of dependence structure on the kernel density
estimation was thoroughly studied by Claeskens and Hall (2002).

The diffusion function can also be consistently estimated when
$\Delta$ is fixed.  In pricing the derivatives of interest rates,
A\"it-Sahalia (1996a) assumes $\mu(x) = k(\alpha-x)$. Using the
kernel density estimator $\hat{f}$ and estimated $\kappa$ and
$\alpha$ from a least-squares method, he applied (\ref{b11}) to
estimate $\sigma(\cdot)$: $\hat{\sigma}^2(x) = 2 \int_0^x
\hat{\mu}(u) \hat{f}(u) du / \hat{f}(x)$. He further established
the asymptotic normality of such an estimator. Gao and King (2004)
propose tests of diffusion models based on the discrepancy between
the parametric and nonparametric estimates of the invariant
density.

The A\"it-Sahalia (1996a) method is a simple one to illustrate
that the volatility function can be consistently estimated for
fixed $\Delta$.  However, we do not expect that it is efficient.
Indeed, we use only the marginal information of the data.  As
shown in (\ref{c1}), almost all information is contained in the
transition density $p_\Delta( \cdot | \cdot )$.  The transition
density can be estimated as in \S 4.2 below, no matter $\Delta$ is
small or large. Since the transition density and drift and
volatility are one-to-one correspondence for the diffusion process
(\ref{b5}). Hence, the drift and diffusion functions can be
consistently estimated via inverting the relationship between the
transition density and the drift and diffusion functions.

There is no simple formula for expressing the drift and diffusion
in terms of the transition density.  The inversion is frequently
carried out via a spectral analysis of the operator $H_\Delta =
\exp(\Delta L)$, where the infinitesimal operator $L$ is defined
as
\[
  L g(x) = \frac{\sigma^2(x)}{2} g''(x) + \mu(x) g'(x).
\]
It has the property:
\[
  L g(x) =  \lim_{\Delta \to 0} \Delta^{-1}
       [ E \{g(X_{t+\Delta})|X_t=x\} - g(x)],
\]
by It\^o's formula (\ref{b10}).  The operator $H_\Delta$ is the
transition operator in that [see also (\ref{b12})]
\[
   H_\Delta g(x) = E \{ g(X_{\Delta}) | X_0 = x \}.
\]
The works of Hansen and Scheinkman (1995), Hansen \etal (1998) and
Kessler and S{\o}rensen (1999) consist of the following idea. The
first step is to estimate the transition operator $H_\Delta$ from
the data. From the transition operator, one can identify the
infinitesimal operator $L$ and hence the functions $\mu(\cdot)$
and $\sigma(\cdot)$. More precisely, let $\lambda_1$ be the
largest negative eigenvalue of the operator $L$ with eigen
function $\xi_1(x)$.  Then, $L \xi_1 = \lambda_1 \xi_1$, or
equivalently $ \sigma^2 \xi_1'' + 2\mu \xi_1' = 2 \lambda_1
\xi_1$.  This gives one equation of $\mu$ and $\sigma$.  Another
equation can be obtained via (\ref{b11}): $(\sigma^2 f)' - 2 \mu f
= 0$.  Solving these two equations we obtain
$$
  \sigma^2(x) = 2\lambda_1 \int_0^x \xi_1(y) f(y) dy / [f(x) \xi_1(x)].
$$
and another explicit expression for $\mu(x)$.
% The drift function $\mu$ can be obtained from
%$$
%\mu(x) = \{\lambda_1 \xi_1(x) - \sigma^2(x)
%\xi_1''(x)/2\}/\xi_1'(x)
%$$
Using the semigroup theory (Theorem IV.3.7, Engel and Nagel 2000),
$\xi_1$ is also an eigen function of $H_\Delta$ with eigenvalue
$\exp(\Delta \lambda_1)$. %See also Chapter 4 of Bass (1997).
Hence, the proposal is to estimate the invariant density $f$ and
the transition density $p_\Delta(y|x)$, which implies the value of
$\lambda_1$ and $\xi_1$.  Gobel \etal (2002) derive the optimal
rate of convergence for such a scheme, using a wavelet basis.  In
particular, they show that for fixed $\Delta$, the optimal rates
of convergence for $\mu$ and $\sigma$ are of orders
$O(n^{-s/(2s+5)})$ and $O(n^{-s/(2s+3)})$, respectively, where $s$
is the degree of smoothness of $\mu$ and $\sigma$.

\subsection{Time-dependent model}

The time dependent model (\ref{b8}) was introduced to accommodate
the possibility of economic changes over time.  The coefficient
functions in (\ref{b8}) are assumed to be slow time-varying and
smooth. Nonparametric techniques can be applied to estimate these
coefficient functions.  The basic idea is to localizing in time,
resulting in a time-domain smoothing.

We first estimate the coefficient functions $\alpha_0(t)$ and
$\alpha_1(t)$.  For each given time $t_0$, approximate the
coefficient functions locally by constants:  $\alpha(t) \approx a$
and $\beta(t) = b$ for $t$ in a neighborhood of $t_0$.  Using the
Euler approximation (\ref{b18}), we run a local regression:
Minimize
\begin{equation}
    \sum \limits_{i=0}^{n-1} (Y_{i \Delta} - a - b X_{i\Delta} )^2
     K_h (i\Delta - t_0) \label{c7}
\end{equation}
with respect to $a$ and $b$.  This results in an estimate
$\hat{\alpha}_0 (t_0) = \hat{a}$ and $\hat{\alpha}_1(t_0) =
\hat{b}$, where $\hat{a}$ and $\hat{b}$ are the minimizer of the
local regression (\ref{c7}).  Fan \etal (2003) suggest using a
one-sided kernel such as $K(u) = (1-u^2)I(-1 < u < 0)$ so that
only the historical data in the time interval $(t_0 - h, t_0)$ are
used in the above local regression. This facilitates forecasting
and bandwidth selection. Our experience shows that there are no
significant differences between nonparametric fitting with
one-sided and two-sided kernels. We opt for local constant
approximations instead of local linear approximations for
estimating time varying functions, since the local linear fit can
create artificial albeit insignificant linear trends when the
underlying functions $\alpha_0(t)$ and $\alpha_1(t)$ are indeed
time-independent. To appreciate this, for constant functions
$\alpha_1$ and $\alpha_2$, a large bandwidth will be chosen to
reduce the variance in the estimation.  This is in essence to fit
a global linear regression for (\ref{c7}).  If the local linear
approximations are used, since no variable selection procedures
have been incorporated in the local fitting (\ref{c7}), the slopes
of the local linear approximations will not be estimated as zero
and hence artificial linear trends will be created for the
estimated coefficients.

The coefficient functions in the volatility can be estimated by
the local approximated likelihood method.  Let
\[
\widehat E_t= \Delta^{-1/2} \{X_{t+\Delta} - X_{t} -
    (\widehat \alpha_0(t) +\widehat \alpha_1(t) X_t)\Delta \}
\]
be the normalized residuals.  Then,
\begin{equation}
 \widehat E_t \approx \beta_0(t) X_t^{\beta_1(t)} \varepsilon_t.
 \label{c8}
\end{equation}
The conditional log-likelihood of $\widehat E_t$ given $X_t$ can
easily be obtained by the approximation (\ref{c8}). Using local
constant approximations and incorporating the kernel weight, we
obtain the local approximated likelihood at each time point and an
estimate of the functions $\beta_0(\cdot)$ and $\beta_1(\cdot)$ at
that time point.
%\begin{equation}
%\ell(\beta_0, \beta_1; t_0) = - \sum_{i=0}^{n-1} K_h(i \Delta - t_0)
%\left ( \log ( \beta_0^2 X_{i\Delta}^{2\beta_1} )
%    +  \frac{\hat{E}_{i\Delta}^2}{\beta_0^2 X_{i\Delta}^{2\beta_1} } \right ).
%\label{c5}
%\end{equation}
%Maximizing (\ref{c5}) with respect to the local parameters $\beta_0$
%and $\beta_1$, we obtain the estimates $\hat{\beta}_0(t_0) =
% \hat{\beta}_0$  and  $\hat{\beta}_1(t_0) = \hat{\beta}_1$.
This type of the local approximated-likelihood method is related
to the generalized method of moments of Hansen (1982), and the
ideas of Florens-Zmirou (1993) and Genon-Catalot and Jacod (1993).

Since the coefficient functions in both return and volatility
functions are estimated using only historical data, their
bandwidths can be selected based on a form of the average
prediction error.  See Fan \etal (2003) for details.  The local
least-squares regression can also be applied to estimate the
coefficient functions $\beta_0(t)$ and $\beta_1(t)$ via the
transformed model [see (\ref{c8})]
\[
\log(\widehat{E}_t^2) \approx 2 \log \beta_0(t) + \beta_1(t) \log(X_t^2)
      + \log(\varepsilon_t^2),
\]
but we do not pursue along this direction, since the local least-squares
estimate is known to be inefficient in the likelihood context and the
exponentiation of an estimated coefficient function of
$\log \beta_0(t)$ is unstable.

A question arises naturally if the coefficients in model
(\ref{b8}) are really time-varying.  This amounts for example to
testing $H_0:  \beta_0(t) = \beta_0$ and $\beta_1(t) = \beta_1$.
Based on the GLR technique, Fan \etal (2003) proposed a formal
test for this kind of problems.

The coefficient functions in the semiparametric model (\ref{b9})
can also be estimated by using the profile approximated-likelihood
method. For each given $\beta_1$, one can estimate easily
$\beta_0(\cdot)$ via the approximation (\ref{c8}), resulting in an
estimate $\hat{\beta}_0(\cdot; \beta_1)$.
%\begin{equation}
% \hat{\beta}_0^2 (t_0; \beta_1) =  \sum_{i=0}^{n-1} K_h(i\Delta - t_0)
% \hat{E}_{i\Delta}^2 |X_{i\Delta}|^{-2\beta_1} \bigg/ \sum_{i=0}^{n-1}
%  K_h(i\Delta - t_0).   \label{c5}
%\end{equation}
Regarding the nonparametric function $\beta_0(\cdot)$ as being
parameterized by $\hat{\beta}_0 (\cdot; \beta_1)$, model
(\ref{c8})  with $\beta_1(t) \equiv \beta_1$ becomes a
``synthesized'' parametric model with unknown $\beta_1$. The
parameter $\beta_1$ can be estimated by the maximum (approximated)
likelihood method.
% Now the coefficient $\beta_1$ can be obtained by maximizing the
%pseudo-likelihood [compare with (\ref{c5})]
%\[
%\ell(\beta_1) = - \sum_{i=0}^{n-1}
%\left (
%  \log \{ \hat\beta_0^2(i\Delta; \beta_1) X_{i\Delta}^{2\beta_1} \}
%      +  \frac{\hat{E}_{i\Delta}^2}{\hat\beta_0^2(i\Delta; \beta_1)
%      X_{i\Delta}^{2\beta_1} } \right ).
%\]
Note that ${\beta}_1$ is estimated by using all the data points,
while $\hat{\beta}_0(t) = \hat\beta_0(t; \hat{\beta}_1)$ is
obtained by using only the local data points. See Fan \etal (2003)
for details.

For other nonparametric methods of estimating volatility in time
inhomogeneous models, see H\"ardle \etal (2003) and Mercurio and
Spokoiny (2003).  Their methods are based on model (\ref{b8}) with
$\alpha_1(t) = \beta_1(t) = 0$.

\subsection{State-domain versus Time-domain smoothing}

So far, we have introduced both state- and time-domain smoothing.
The former relies on the structural invariability implied by the
stationarity assumption and uses pre-dominantly on the (remote)
historical data.  The latter uses the continuity of underlying
parameters and concentrates basically on the recent data.  This
can be illustrated in Figure 7, using the yields of the 3-month
Treasury bills from January 8, 1954 to July 16, 2004, sampled at
weekly frequency. On December 28, 1990, the interest rate is about
6.48\%.  To estimate the drift and diffusion around $x=6.48$, the
state-domain focuses on the dynamics where interest rates are
around 6.48\%, the horizontal bar with interest rates falling in
$6.48\% \pm .25\%$. The estimated volatility is basically the
sample standard deviation of the differences $\{X_{i\Delta} -
X_{(i-1)\Delta}\}$ within this horizontal bar. On the other hand,
the time-domain smoothing focuses predominately on the recent
history, say one year, as illustrated in the figure. The
time-domain estimate of volatility is basically a sample standard
deviation within the vertical bar.

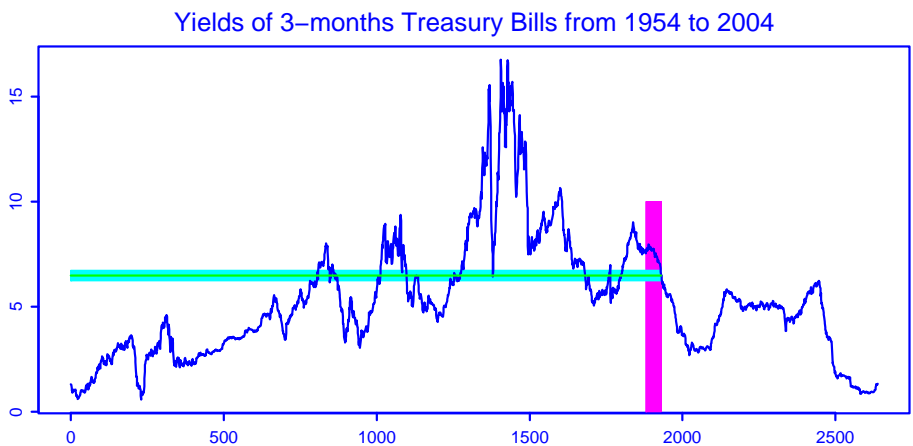
\fig{figs}{Fig7}{Illustration of the time and state-domain
smoothing using the yields of 3-month Treasury bills.  The
state-domain smoothing localizing in the horizontal bars, while
the time-domain smoothing concentrating in the vertical bars.}

For a given time series, it is hard to say which estimate is
better.  This depends on the underlying stochastic processes and
also on the time when the forecast to be made.  If the underlying
process is continuous and stationary such as model (\ref{b5}),
both methods are applicable.  For example, standing on December
28, 1990, one can forecast the volatility by using the sample
standard deviation in either the horizontal bar or vertical bar.
However, the estimated precision depends on the local data. Since
the sample variance is basically linear in the squared
differences, the standard errors of both estimates can be assessed
and used to guide the forecasting.
%The standard errors depend ont he time of forecasting.

For stationary diffusion processes, it is possible to integrate
both the time-domain and state-domain estimates.  Note that the
historical data (with interest rates in $6.48\% \pm .25\%$) are
far apart in time (except the last piece, which can be ignored in
the state-domain fitting) from the data used in the time-domain
smoothing (vertical bar).  Hence, these two estimates are nearly
independent.  The integrated estimate is a linear combination of
these two nearly independent estimates.  The weights can easily be
chosen to minimize the variance of the integrated estimator, by
using the assessed standard errors of the state- and time-domain
estimators. This forms a dynamically integrated predictor for
volatility estimation, as the optimal weights change over time.
%The idea can also be applied to forecast the return function.

\subsection{Continuously observed data}

At theoretical level, one may also examine the problem of
estimating the drift and diffusion functions assuming the whole
process is observable up to time $T$.  Let us assume again that
the observed process $\{X_t\}$ follows SDE (\ref{b5}). In this
case, $\sigma^2(X_t)$ is the derivative of the quadratic variation
process of $X_t$ and hence is known up to time $T$. By
(\ref{b11}), estimating the drift function $\mu(x)$ is equivalent
to estimating the invariant density $f$.  In fact,
\begin{equation}
  \mu(x) = [\sigma^2(x) f(x) ]'/[2f(x)]. \label{c9}
\end{equation}

The invariant density $f$ can easily be estimated by the kernel
density estimation.  When $\Delta \to 0$, the summation in
(\ref{c6}) converges to
\begin{equation}
    \hat{f}(x) = T^{-1} \int_0^T K_h (X_t - x) dt. \label{c10}
\end{equation}
This forms a kernel density estimate of the invariant density
based on the continuously observed data.  Thus, an estimator for
$\mu(x)$ can be obtained by substituting $\hat{f}(x)$ into
(\ref{c9}). Such an approach has been employed by Kutoyants (1998)
and Dalalyan and Kutoyants (2000, 2003).  They established sharp
asymptotic minimax risk for estimating the invariant density $f$
and its derivative, as well as the drift function $\mu$. In
particular, the functions $f$, $f'$ and $\mu$ can be estimated
with rate $T^{-1/2}$, $T^{-2s/(2s+1)}$ and $T^{-2s/(2s+1)}$,
respectively, where $s$ is the degree of smoothness of $\mu$.
These are the optimal rates of convergence.

An alternative approach is to estimate the drift function directly
from (\ref{c3}).  By letting $\Delta \to 0$, one can easily obtain
a local linear regression estimator for continuously observed
data, which admits a similar form to (\ref{c3}) and (\ref{c10}).
This is the approach that Spokoiny (2000) used. He
% applies the local linear
%estimator to estimate the drift function. By taking the limit as
%$\Delta \to 0$, the estimator in (\ref{c3}) converges to
%\begin{equation}
%\hat{\mu}(x) = \frac{ S_2(x) \int_0^T K_h (X_t -x) dX_t  -
%    S_1(x) \int_0^T K_h(X_t
%-x)(X_t -x) dX_t}{S_2(x) S_0(x) - S_1(x)^2}, \label{c8}
%\end{equation}
%where $S_j(x) = \int_0^T K_h (X_t - x) (X_t -x)^j dt$.
showed that this estimator attains the optimal rate of convergence
and established further a data-driven bandwidth such that the
local linear estimator attains adaptive minimax rates.

%The above local linear estimator can be derived directly from the local
%modeling technique (Fan and Gijbels, 1996).  Note that
%\[
%E \left \{ \int_0^T f(X_t)^2 dt - 2 \int_{0}^T f(X_t) dX_t \right \}
%   = E \int_0^T (f(X_t) - \mu(X_t))^2 dt - \int_0^T \mu(X_t)^2 dt
%\]
%is minimized at $f = \mu$.  Thus, approximating $f(X_t)$ locally by
%$f(x) + f'(x) (X_t -x)$ for $X_t$ around $x$, one would minimize
%\begin{equation}
%  \int_0^T \{a + b(X_t -x)\}^2 K_h(X_t -x) dt - 2 \int_{0}^T \{a + b(X_t
%  -x)\} K_h(X_t - x) dX_t  \label{c9}
%\end{equation}
%with respect to $a$ and $b$.  Minimizing (\ref{c9}) results in the same
%estimator as (\ref{c8}).

\section{Estimation of state price densities and transition densities}

State-price density (SPD) is the probability density of the value
of an asset under the risk-neutral world (\ref{b14}) [see Cox and
Ross (1976)] or equivalent martingale measure (Harrison and Kreps,
1979).  It is directly related to the pricing of financial
derivatives.  It is the transition density of $X_T$ given $X_0$
under the equivalent martingale $Q$.  The SPD does not depend on
the payoff function and hence it can be used to evaluate other
illiquid derivatives, once it is estimated from more liquid
derivatives.  On the other hand, the transition density
characterizes the probability law of a Markovian process and hence
is useful for validating Markovian properties and parametric
models.

\subsection{Estimation of state price density}

For some specific models, the state price density can be formed
explicitly.  For example, for the GBM (\ref{b1}) with a constant
risk-free rate $r$, according to (\ref{b16}), the SPD is
log-normal, with mean $\log x_0 + (r - \sigma^2)/(2T)$ and
variance $\sigma^2$.

Assume that the SPD $f^*$ exists.  Then, the European call option can be
expressed as
\[
   C = \exp(-\int_0^T r_s ds ) \int_K^\infty (x - K) f^*(x) dx.
\]
See (\ref{b14}) (we have changed the notation from $P_0$ to $C$ to
emphasize the price of the European call option). Hence,
\begin{equation}
   f^*(K) = \exp(\int_0^T r_s ds ) \frac{\partial^2 C}{\partial K^2}.
   \label{d1}
\end{equation}
This was observed by Breeden and Litzenberger (1978).  Thus, the
state price density can be estimated from the European call
options with different strike prices.  With the estimated state
price density, one can price new or less-liquid securities such as
over the counter derivatives or nontraded options, using formula
(\ref{b14}).

In general, the price of an European call option depends on the
current stock price $S$, the strike price $K$, the time to
maturity $T$, the risk-free interest rate $r$ and dividend yield
rate $\delta$.  It can be written as $C(S, K, T, r, \delta)$.  The
exact form of $C$, in general, is hard to determine, unless we
assume the Black-Scholes model. Based on historical data $\{(C_i,
S_i, K_i, T_i, r_i, \delta_i), i = 1, \cdots, n\}$, where $C_i$ is
the $i^{th}$ traded-option price with associated characteristics
$(S_i, K_i, T_i, r_i, \delta_i)$, A\"it-Sahalia and Lo (1998) fit
the following nonparametric regression
$$
   C_i = C(S_i, K_i, T_i, r_i, \delta_i) + \varepsilon_i
$$
to obtain an estimate of the function $C$ and hence the SPD $f^*$.

Due to the curse of dimensionality, the five dimensional
nonparametric function can not be estimated well with practical
range of sample sizes. A\"it-Sahalia and Lo (1998) realized that
and proposed a few dimensionality reduction methods.  First, by
assuming that the option price depends only on the futures price
$F = S \exp( (r - \delta)T)$, namely,
\[
    C(S, K, T, r, \delta) = C(F, K, T, r)
\]
(the Black-Scholes formula satisfies such an assumption), they
reduced the dimensionality from 5 to 4.  By assuming further that
the option-pricing function is homogeneous of degree one in $F$
and $K$, namely,
\[
  C(S, K, T, r, \delta) = K C(F/K, T, r),
\]
they reduced the dimensionality to 3. A\"it-Sahalia and Lo (1998)
imposed a semiparametric form on the pricing formula:
\[
   C(S, K, T, r, \delta) = C_{\mbox{\scriptsize BS}} (F, K, T, r,
   \sigma(F, K, T)),
\]
where $C_{\mbox{\scriptsize BS}} (F, K, T, r, \sigma)$ is the
Black-Scholes pricing formula given in (\ref{b17}) and $\sigma(F,
K, T)$ is the implied volatility, computed by inverting the
Black-Scholes formula.  Thus, the problem becomes
nonparametrically estimating the implied volatility function
$\sigma(F, K, T)$.  This is estimated by using a nonparametric
regression technique from historical data, namely
$$
  \sigma_i = \sigma(F_i, K_i, T_i) + \varepsilon_i,
$$
where $\sigma_i$ is the implied volatility of $C_i$ by inverting the
Black-Scholes formula.
By assuming further that $\sigma(F, K, T) = \sigma(F/K, T)$, the
dimensionality is reduced to 2.
This is one of the options in A\"it-Sahalia and
Lo (1998).

The state price density $f^*$ is non-negative and hence the
function $C$ should be convex in the strike price $K$.
A\"it-Sahalia and Duarte (2003) propose to estimate the option
price under the convexity constraint, using a local linear
estimator. See also H\"ardle and Yatchew (2002) for a related
approach.

\subsection{Estimation of transition densities}

The transition density of a Markov process characterizes the law
of the process, except the initial distribution.  It provides
useful tools for checking whether or not such a process follows a
certain SDE and for statistical estimation and inferences.  It is
the state price density of the price process under the risk
neutral world. If such a process were observable, the state price
density can be estimated using the methods to be introduced.

Assume that we have a sample $\{X_{i\Delta}, i = 0, \cdots, n\}$
from model (\ref{b5}).   The ``double-kernel'' method of Fan \etal
(1996) is to observe that
\begin{equation}
   E \{ W_{h_2} (X_{i\Delta} - y)|X_{(i-1)\Delta}=x\} \approx p_\Delta(y|x),
   \quad \mbox{as $h_2 \to 0$},   \label{d2}
\end{equation}
for a kernel function $W$. Thus, the transition density
$p_\Delta(y|x)$ can be regarded approximately as the nonparametric
regression function of the response variable $W_{h_2} (X_{i\Delta}
- y)$ on $X_{(i-1)\Delta}$.  An application of the local linear
estimator (\ref{c3}) yields
\begin{equation}
\hat{p}_\Delta (y|x) = \sum_{i=1}^n K_n(X_{(i-1)\Delta} - x, x)
W_{h_2} (X_{i\Delta} - y), \label{d3}
\end{equation}
where the equivalent kernel $K_n(u, x)$ was defined in (\ref{c4}).
Fan \etal (1996) establish the asymptotic normality of such an
estimator under stationarity and $\rho$-mixing conditions
[necessarily decaying at geometric rate for SDE (\ref{b5})], which
gives explicitly the asymptotic bias and variance of the
estimator.  See also \S6.5 of Fan and Yao (2003). The
cross-validation idea of Rudemo (1982) and Bowman (1984) can be
extended to select bandwidths for estimating conditional
densities. See Fan and Yim (2004) and Hall \etal (2004).

The transition distribution can be estimated by integrating
the estimator (\ref{d3}) over $y$.  Alternative estimators can be
obtained by an application of the local logistic regression and adjusted
Nadaraya-Watson method of Hall \etal (1999).

Early references on the estimation of the transition distributions and
densities include Roussas (1967, 1969) and Rosenblatt (1970).

\subsection{Inferences based on transition densities}

With the estimated transition density, one can now verify whether
parametric models such as (\ref{b1})--(\ref{b3}), (\ref{b6}) are
consistent with the observed data. Let $p_{\Delta, \theta}(y|x)$
be the transition density under a parametric diffusion model. For
example, for the CIR model (\ref{b7}), the parameter $\theta =
(\kappa, \alpha, \sigma)$.  As in (\ref{c1}), ignoring the initial
value $X_{0}$, the parameter $\theta$ can be estimated by
maximizing
$$
  \ell(p_{\Delta, \theta}) =
  \sum_{i=1}^n \log p_{\Delta, \theta} (X_{i\Delta}|X_{(i-1) \Delta}).
$$
Let $\hat{\theta}$ be the maximum likelihood estimator. By the
spirit of the GLR of Fan \etal (2001), the GLR test for the null
hypothesis $ H_0:  p_\Delta(y|x) = p_{\Delta, \theta}(y|x)$ is
\[
    \mbox{GLR} =  \ell(\hat{p}_{\Delta}) - \ell(p_{\Delta, \hat \theta}),
\]
where $\hat{p}$ is a nonparametric estimate of the transition
density.  Since the transition density can not be estimated well over the
region where data are sparse (usually at boundaries of the process), we
need to truncate the nonparametric (and simultaneously parametric)
evaluation of the likelihood at appropriate intervals.

In addition to employing the GLR test, one can also compare
directly the difference between the parametric and nonparametric
fits, resulting in test statistics such as $ \| \hat{p}_\Delta  -
p_{\Delta, \hat \theta} \|^2$ and $ \| \hat{P}_\Delta  -
P_{\Delta, \hat \theta} \|^2$ for an appropriate norm $\|\cdot
\|$, where $\hat{P}_\Delta$ and $P_{\Delta, \hat \theta}$ are the
estimates of the cumulative transition distributions under
respectively the parametric and nonparametric model. An
alternative method is to apply the GLR of Fan \etal (2001) to
separately test the forms of the drift and diffusion, as in Fan
and Zhang (2003). The transition density approach appears more
elegant as it checks simultaneously the forms of drift and
diffusion, but more computationally intensive.  In comparisons
with the invariant density-based approach of Arapis and Gao
(2003), it is consistent against a much larger family of
alternatives.

One can also use the transition density to test whether an observed
series is Markovian (from personal communication with Yacine
A\"it-Sahalia).  For example, if a process
$\{X_{i \Delta}\}$ is Markovian, then
\[
   p_{2 \Delta}(y|x) = \inte p_{\Delta}(y|z) p_{\Delta}(z|x) dz.
\]
Thus, one can use the distance between $\hat{p}_{2 \Delta}(y|x)$
and $\inte \hat{p}_{\Delta}(y|z) \hat{p}_{\Delta}(z|x) dz$ as a
test statistic.

Transition density can also be used for parameter estimation.  One
possible approach is to find the parameter to minimize the
distance $\| \hat{P}_\Delta  - P_{\Delta, \theta}\|$.  In this
case, the bandwidth should be chosen to optimize the performance
for estimating $\theta$.

\section{Concluding remark}

Enormous efforts in financial econometrics have been made in
modeling the dynamics of stock prices and bond yields.  There are
directly related to pricing derivative securities, proprietary
trading and portfolio management. Various parametric models have
been proposed to facilitate mathematical derivations.  They have
risks that misspecifications of models lead to erroneous pricing
and hedging strategies.  Nonparametric models provide a powerful
and flexible treatment.  They aim at reducing modeling biases by
increasing somewhat the estimation variances.   They provide an
elegant method for validating or suggesting a family of parametric
models.

The versatility of nonparametric techniques in financial
econometrics has been demonstrated in this paper.  They are
applicable to various aspects of diffusion models:  drift,
diffusion, transition densities, and even state price densities.
They allow us to examine whether the stochastic dynamics for
stocks and bonds are time varying and whether famous parametric
models are consistent with empirical financial data.  They permit
us to price illiquid or non-traded derivatives from liquid
derivatives.

The applications of nonparametric techniques in financial
econometrics are far wider than what has been presented.  There
are several areas where nonparametric methods have played a
pivotal role.  One example is to test various versions of capital
asset pricing models (CAPM) and their related stochastic discount
models (Cochrane, 2001).  See for example the research manuscript
by Chen and Ludvigson (2003) in this direction. Another important
class of models are stochastic volatility models (Barndoff-Neilsen
and Shephard, 2001 and Shephard 2004), where nonparametric methods
can be applied.  The nonparametric techniques have been
prominently featured in the RiskMetrics of J. P. Morgan.  It can
be employed to forecast the risks of portfolios.  See, for
example, Jorion (2000), A\"i-Sahalia and Lo (2000), Chen and Tang
(2003), Fan and Gu (2003) and Chen (2004).

\begin{singlespace}

\end{singlespace}

\end{document}